\documentclass[aap,MSNbibl,nameyear,dvips]{arximspdf}
\usepackage{constants}

\doi{10.1214/12-AAP848} %kopijuoti is PTS
\volume{23}
\issue{2}
\pubyear{2013}
\firstpage{617}
\lastpage{636}

\makeatletter

\newtheorem{theorem}{Theorem}[section]

\newtheorem{lemma}{Lemma}[section]
\newproclaim{definition}{Definition}[section]
\newproclaim{example}{Example}[section]
\newproclaim{remark}{Remark}[section]
\newcommand{\eqref}[1]{(\ref{#1})}
\renewcommand{\citeyearpar}[1]{(\citeyear{#1})}
\makeatother

\begin{document}
\begin{frontmatter}

\title{A Berry--Esseen bound with applications to vertex degree counts
in the Erd\H{o}s--R\'enyi random graph}
\runtitle{Counts in the Erd\H{o}s--R\'enyi random graph}

\begin{aug}
\author[A]{\fnms{Larry} \snm{Goldstein}\corref{}\thanksref{t1}\ead[label=e1]{larry@math.usc.edu}}
\thankstext{t1}{Supported in part by NSA Grant H98230-11-1-0162.}
\runauthor{L. Goldstein}
\affiliation{University of Southern California}
\address[A]{Department of Mathematics\\
University of Southern California\\
Kaprielian Hall, Room 108\\
3620 Vermont Avenue\\
Los Angeles, California 90089-2532\\
USA\\
\printead{e1}} %adresu isvedimo komanda gale!
\end{aug}

% HISTORY:
\received{\smonth{5} \syear{2010}}
\revised{\smonth{12} \syear{2011}}

% ABSTRACT
%
\begin{abstract}
Applying Stein's method, an inductive technique and size bias coupling
yields a Berry--Esseen theorem for normal approximation without the
usual restriction that the coupling be bounded. The theorem is applied
to counting the number of vertices in the Erd\H{o}s--R\'enyi random
graph of a given degree.
\end{abstract}

% KEYWORDS
%

\end{frontmatter}
%

%s1 #&#
\section{Introduction}

We present a new Berry--Esseen theorem for sums $Y$ of dependent
variables by combining Stein's method, size bias couplings and the
inductive technique of Bolthausen \citeyearpar{Bolthausen84} originally
developed for the combinatorial central limit theorem. We apply the
theorem to asses the accuracy of the normal approximation to the
distribution of the number of vertices of degree $d$ in the classical
Erd\H{o}s--R\'enyi \citeyearpar{Erdos59} random graph $G_n$ having $n$
vertices connected by independent edges with common success probability
depending on $n$ and a parameter $\theta$. Over the range of parameters
considered, the theorem yields a bound that is the same up to constants
as the one obtained earlier by \citet{Barbour89} for the weaker smooth
function metric (\ref{defd1}).

Stein's method [\citeauthor{Stein72} (\citeyear{Stein72,Stein86})] often proceeds by
coupling a random variable $Y$ of interest to a related variable $Y'$,
using, for example, the method of exchangeable pairs, size bias
couplings or zero bias couplings; for an overview see \citet{Chen10}.
The chief innovation here is the removal of an inconvenient restriction
present in a number of results that provide Kolmogorov distance bounds
using Stein's method, that the difference $|Y-Y'|$ between $Y$ and the
coupled $Y'$ be bounded almost surely by a constant. Through the use of
an unbounded coupling, in Theorem~\ref{thmdegree} we are able to
extend the previous work by \citet{Kordecki90} on the number of
isolated, or degree zero, vertices of $G_n$ to all positive
degrees.\vadjust{\goodbreak}

To describe Theorem~\ref{thmmain}, our general result, recall that for
a nonnegative random variable $Y$ with finite, nonzero mean $\mu$, we
say that $Y^s$ has the $Y$-size bias distribution if
%
%e1 #&#
\begin{equation}\label{introszb}
E[Yf(Y)]=\mu E[f(Y^s)]
\end{equation}
for all functions $f$ for which
these expectations exist.

In employing the size bias version of Stein's method [see \citet
{Baldi89}, \mbox{\citet{Goldstein96}} and \citet{Chen10}], the goal is to
construct, on the same space as $Y$, a variable $Y^s$ with the $Y$-size
bias distribution such that $Y$ and $Y^s$ are close is some sense.
Previous applications of the size bias coupling technique for obtaining
Berry--Esseen bounds by Stein's method, requiring that $|Y^s-Y|$ be
bounded, include \citet{Go05}, \citet{Goldstein10} and \citet{GoZh}.

Let $\mathbb{N}=\{1,2,\ldots\}$ and $\mathbb{N}_0=\mathbb{N} \cup\{
0\}
$. Our abstract framework consists of random elements indexed by $n \ge
n_0$ for some $n_0 \in\mathbb{N}_0$ whose distributions $\mathcal{
L}_\theta(\cdot)$ depend on $n$, left implicit when clear from context,
and a parameter $\theta$ in a topological space $\Theta_n$. We also
assume that $\Theta_n$ is endowed with a $\sigma$-algebra, taken to be
the one generated by the collection of open sets unless specified otherwise.

In our application the parameter $\theta$ lies in a subset $\Theta_n$
of the real numbers $\mathbb{R}$ and interest centers on the
distributions of the nonnegative random variables $Y_n$ counting the
number of degree $d \in\mathbb{N}$ vertices of the Erd\H{o}s--R\'enyi
random graph~$G_n$.
%For notational ease we may drop the subscript $\theta$ when it is not
%required for emphasis.
For sums of exchangeable indicator variables such as $Y_n$, Lemma \ref
{sblem} below says, essentially, that to construct a variable $Y_n^s$
with the $Y_n$-size bias distribution, one chooses an indicator
uniformly and sets it to one if it was not so already, and then
``adjusts'' the remaining indicators, if necessary, to have their
original distribution given that the selected indicator is one.
Applying Lemma~\ref{sblem} when $Y_n$ counts the number of vertices in
$G_n$ having degree $d$ results in the construction of \citet
{Barbour92}, where nothing is changed if a uniformly chosen vertex
already has degree $d$, and otherwise edges to the chosen vertex are
added if the vertex has degree less than $d$, or removed if it has
degree in excess of $d$. As it is possible that the chosen vertex has,
say, $n-1$ edges, the resulting coupling fails to be bounded in $n$.
Nevertheless, when there is only a small probability that a very large
number of edges will need to be added or removed, the coupling can be
controlled using moments on bounds $K_n$ that satisfy $|Y_n^s-Y_n| \le K_n$.

After coupling, the second ingredient in our method has an inductive
flavor. We construct a variable $V_n$ such that its distribution,
conditional on a collection $J_n$ of random elements, is that of $Y_n$
reduced in size by some ``small'' amount $L_n$, with parameter $\psi
_{n,\theta}$ ``close'' to the original $\theta$. Formally, we require that
%
%e2 #&#
\begin{equation}\label{introYncond}
\mathcal{ L}_\theta(V_n|J_n)=\mathcal{ L}_{\psi_{n,\theta}}(Y_{n-L_n})
\end{equation}
hold on an event where the size of $L_n$ is controlled, and that a
bound $B_n$ on the absolute difference $|Y_n-V_n|$ not be ``too large.''
As bounds to the normal for $Y_n$ can be expressed in terms of
quantities that include bounds to the normal for reduced versions of
the same problem, a recursive inequality for the sought after bound can
be produced.

In the graph degree problem, $V_n$ counts the number of degree $d$
vertices in the graph obtained by removing a uniformly chosen vertex
from $G_n$, along with all its incident edges, and the set $J_n$
consists of the identity of the chosen vertex, and its degree.
Conditionally on $J_n$, the graph that remains is an Erd\H{o}s--R\'enyi
graph on the reduced vertex set, with the same connectivity as before.
As with the bound $K_n$, it is not required that $B_n$ be almost surely
bounded by a constant; though $|Y_n-V_n|$ may be large in the graph
degree problem, it is unlikely that it will be.

Tension exists in choosing the set $J_n$ that appears in the
conditioning equality~(\ref{introYncond}). In order to reduce the
larger problem to a smaller one so that induction may be applied,
working conditionally we must be able to treat the bounds $K_n$ and~$B_n$,
and the parameters of the reduced problem, $L_n$ and $\psi
_{n,\theta}$, as constants. Hence we require that these variables be
measurable with respect to $\mathcal{ F}_n$, the \mbox{$\sigma$-algebra} generated
by the conditioning collection $J_n$. Though this restriction
necessitates that $\mathcal{ F}_n$ be large enough to contain, say,
information on $Y_n^s-Y_n$, it must also be small enough so that $L_n$
and $B_n$ are not too large, and that the conditioning ``leaves enough
randomness'' to yield a useful recursion for the ultimate bound.

At the heart of our main result, and Stein's method for normal
approximation, is the characterization that $Z$ is a standard normal
random variable if and only if
\[
E[Zf(Z)]=E[f'(Z)]
\]
for all absolutely continuous functions $f$ for which the above
expectations exist. This characterization leads to the Stein equation,
when, given a test function $h$ on which to evaluate the difference
$Eh(W)-Eh(Z)$ between the expectation of the random variable $W$ of
interest and the standard normal $Z$, one solves
\[%\label{introsteq}
f'(w)-wf(w)=h(w)-Eh(Z)
\]
for $f$. Using $f$, one evaluates this difference by substituting $W$
for $w$, and takes expectation on the left-hand side, rather than the
right. Though we focus on manipulation of the Stein equation using the
size bias coupling, many variations are possible; see \citet{Chen10}
for an overview.\vadjust{\goodbreak}

Throughout, for $n_0 \in\mathbb{N}$ and all $n \ge n_0$ and $\theta
\in\Theta_n$, we let $\mu_{n,\theta}=E_\theta Y_n$ and $\sigma
_{n,\theta}^2=\operatorname{Var}_\theta(Y_n)$ indicate the
mean and variance of $Y_n$ under $\mathcal{ L}_\theta$. The value
$r_{n,\theta}$ appearing in Theorem~\ref{thmmain}
is a function that determines the quality of the bound to the normal,
while the sequence $s_{n,\theta}$ is used to control $L_n$, and hence
the size of the smaller subproblem $V_n$ related to $Y_n$. Without
further mention, $\mu_{n,\theta},\sigma_{n,\theta}^2$ and
$r_{n,\theta
}$ are assumed to be measurable in $\theta\in\Theta_n$, a condition
satisfied for all natural examples, including the one considered here.
To avoid repetition, the distribution of random variables indicated
after $\theta\in\Theta_n$ has been fixed is with respect to $\mathcal{
L}_\theta$. The random variable $Z$ will always denote the standard normal.
% so that quantities such as $\mu_{n-L_n,\psi_{n,\theta}$ are $\mathcal{
%F}_n$ measurable random variables.

To familiarize the reader with the conditions of Theorem \ref
{thmmain}, toward the end of this section we present its application
in the simple case where a bounded size bias coupling of $Y_n^s$ to
$Y_n$ exists.\vspace*{-1pt}

%th1.1 #&#
\begin{theorem}
\label{thmmain}
%For some $n_0 \in\mathbb{N}_0$ and all $n \ge n_0$
%and $\theta\in\Theta_n$, let
%$Y_n$ be a nonnegative random variable with mean $\mu_{n,\theta}$ and
%positive variance $\sigma_{n,\theta}^2$, and set
For some $n_0 \in \mathbb{N}_0$ and all $n \ge n_0$, let
$Y_n$ be a nonnegative random variable with mean $\mu_{n,\theta}=E_\theta Y_n$
and positive variance $\sigma_{n,\theta}^2=\operatorname{Var}_\theta(Y_n)$ for all $\theta \in \Theta_n$, and set\vspace*{-1pt}
%
%e3 #&#
\begin{equation}\label{defWntheta}
W_{n,\theta}=\frac{Y_n-\mu_{n,\theta}}{\sigma_{n,\theta}},
\end{equation}
the standardized value of $Y_n$. Let $r_{n,\theta}$ be positive for all
$n \ge n_0$ and all $\theta\in\Theta_n$, and for all $r \ge0$ let
\[%\label{defThetanr}
\Theta_{n,r} = \{\theta\in\Theta_n\dvtx  r_{n,\theta} \ge r\}.
\]
Assume there exists $r_1>0$ and $n_1 \ge n_0$ such that
%
%e4 #&#
\begin{equation}\label{bdawayzerornt}
\max_{n_0 \le n < n_1} \sup_{\theta\in\Theta_{n,r_1}} r_{n,\theta
} <
\infty.
\end{equation}
Further, suppose that for all $n \ge n_1$ and $\theta\in\Theta
_{n,r_1}$, there exist random variables $Y_n^s, K_n, L_n,$ $\psi
_{n,\theta}, V_n$ and $B_n$ on the same space as $Y_n$, and a $\sigma
$-algebra $\mathcal{ F}_n$, generated by a collection of random elements
$J_n$, such that the following conditions hold:

\begin{longlist}[1.]
\item[1.] \hypertarget{Psinisrootn} The random variable $Y_n^s$ has the
$Y_n$-size bias distribution, and
%
%e5 #&#
\begin{eqnarray}\label{defpsin}
&&\Psi_{n,\theta}=\sqrt{\operatorname{Var}_\theta\bigl(E_\theta
(Y_n^s-Y_n|Y_n)\bigr)} \qquad\mbox{satisfies}
\nonumber
\\[-8pt]
\\[-8pt]
\nonumber
&&\qquad\sup_{n \ge n_1,\theta\in
\Theta_{n,r_1}}\frac{r_{n,\theta}\mu_{n,\theta} \Psi_{n,\theta
}}{\sigma
_{n,\theta}^2} < \infty.
\end{eqnarray}

\item[2.] \hypertarget{assumptionKn} The random variable $K_n$ is $\mathcal{
F}_n$-measurable, $|Y_n^s-Y_n| \le K_n$ and
%
%e6 #&#
\begin{equation}\label{assumptionKninequality}
\sup_{n \ge n_1,\theta\in\Theta_{n,r_1}}\frac{r_{n,\theta}\mu
_{n,\theta} E_\theta[( 1+ |W_{n,\theta}|
)K_n^2]}{\sigma
_{n,\theta}^3}<\infty
\end{equation}
with $W_{n,\theta}$ as given in (\ref{defWntheta}).

\item[3.] \hypertarget{assumptioncondlawYnvgivenK} The random variable
$L_n$ takes values in $\{0,1,\ldots,n\}$, there exists a positive
integer valued sequence $\{s_{n,\theta}\}_{n \ge n_1}$ satisfying
$n-s_{n,\theta} \ge n_0$, the variables $L_n$ and $\psi_{n,\theta}$ are
$\mathcal{ F}_n$-measurable,
for some $F_{n,\theta} \in \mathcal{F}_n$ satisfying $F_{n,\theta} \subset \{L_n \le s_{n,\theta}\}$,
%there exists $F_{n,\theta} \in\mathcal{ F}_n,
%F_{n,\theta} \subset\{L_n \le s_{n,\theta}\}$ such that
%
%e7 #&#
\begin{eqnarray}\label{VnconditionedJn}
\psi_{n,\theta} \in\Theta_{n-L_n} \quad\mbox{and}\quad \mathcal{ L}_\theta
(V_n|J_n)=\mathcal{ L}_{\psi_{n,\theta}}(Y_{n-L_n})\qquad \mbox{on
$F_{n,\theta}$}
\end{eqnarray}
and
%
%e8 #&#
\begin{equation}\label{musigKLs}
\sup_{n \ge n_1,\theta\in\Theta_{n,r_1}}\frac{r_{n,\theta}^2\mu
_{n,\theta}}{\sigma_{n,\theta}^3} E_\theta\bigl[ K_n^2 \bigl(1-\mathbf{1}(F_{n,\theta})\bigr)\bigr] < \infty.
\end{equation}

\item[4.] \hypertarget{unifvarbnd}
There exists $\{c_1, c_2\} \subset(0,\infty)$ such that
\[%\label{unifvarbnd12}
\sigma^2_{n,\theta} \le c_1 \sigma_{n-L_n,\psi_{n,\theta}}^2
\quad\mbox{and}\quad
r_{n,\theta} \le c_2 r_{n-L_n,\psi_{n,\theta}}\qquad
\mbox{on $F_{n,\theta}$.}
\]

\item[5.] \hypertarget{assumptionYnmYnVnprime} The random variable $B_n$ is
$\mathcal{ F}_n$-measurable, $|Y_n-V_n| \le B_n$ and
%
%e9 #&#
\begin{equation}\label{inequalityYnmYnVnprime}
\sup_{n \ge n_1,\theta\in\Theta_{n,r_1}}\frac{r_{n,\theta}^2\mu
_{n,\theta}
E_\theta[K_n^2 B_n]
}{\sigma_{n,\theta}^4}<\infty.
\end{equation}

\item[6.] \hypertarget{Thetanr-compact} Either:
\begin{longlist}[(a)]

\item[(a)] there exists $l_{n,0} \in\{0,\ldots,n\}$ such that $P_\theta
(L_n=l_{n,0})=1$ for all $\theta\in\Theta_{n,r_1}$,
or

\item[(b)] the set $\Theta_{n,r_1}$ is a compact subset of $\Theta_n$, and
the functions of $\theta$
%
%e10 #&#
\begin{equation}\label{deftnthetal}
t_{n,\theta,l}= E_\theta\biggl( \frac{K_n^2}{E_\theta K_n^2}\mathbf{1}(L_n=l)\biggr), \qquad\mbox{$l \in\{0,1,\ldots,n\}$}
\end{equation}
are continuous on $\Theta_{n,r_1}$ for $l \in\{0,1,\ldots,s_n\}$ where
$s_n=\sup_{\theta\in\Theta_{n,r_1}}s_{n,\theta}$.

\end{longlist}
\end{longlist}

Then there exists a constant $C$ such that for all $n \ge n_0$ and
$\theta\in\Theta_n$,
%
%e11 #&#
\begin{equation}\label{thmmainstatement}
\sup_{z \in\mathbb{R}}|P_\theta(W_{n,\theta} \le z) - P(Z \le z)|
\le
C/r_{n,\theta}.
\end{equation}
%
%where $Z$ is a standard normal random variable.
\end{theorem}

When higher moments exist a number of the conditions of the theorem
may be verified using standard inequalities. In particular, by the
Cauchy--Schwarz inequality a sufficient condition for (\ref
{assumptionKninequality}) is
%
%e12 #&#
\begin{equation}
\label{kn-sum2-conv}
\sup_{n \ge n_1,\theta\in\Theta_{n,r_1}}\frac{r_{n,\theta}\mu
_{n,\theta} k_{n,\theta,4}^{1/2}}{\sigma_{n,\theta}^3}<\infty\qquad
\mbox{where }   k_{n,\theta,m}=E_\theta K_n^m,
\end{equation}
and, when $F_{n,\theta}=\{L_{n} \le s_{n,\theta}\}$, a
sufficient condition for (\ref{musigKLs}) is
\[%\label{musigKLssufficient}
\sup_{n \ge n_1,\theta\in\Theta_{n,r_1}}\frac{r_{n,\theta}^2\mu
_{n,\theta} k_{n,\theta,4}^{{1}/{2}} l_{n,\theta,2}^{
{1}/{2}}}{\sigma_{n,\theta}^3 s_{n,\theta}} < \infty\qquad\mbox{where }
l_{n,\theta,m}=E_\theta L_n^m,
\]
since, additionally using the Markov inequality yields
\begin{eqnarray*}
E_\theta[ K_n^2 \mathbf{1}(L_n > s_{n,\theta})]
&\le &k_{n,\theta,4}^{1/2}P_\theta(L_n>s_{n,\theta})^{{1}/{2}}
= k_{n,\theta,4}^{1/2}P_\theta(L_n^2>s_{n,\theta}^2)^{{1}/{2}}\\
&\le&\frac{k_{n,\theta,4}^{{1}/{2}}l_{n,\theta,2}^{
{1}/{2}}}{s_{n,\theta}}.
\end{eqnarray*}
Similarly, a sufficient condition for (\ref{inequalityYnmYnVnprime}) is
%
%e13 #&#
\begin{eqnarray}\label{inequalityYnmYnVnprimesufficient}
\sup_{n \ge n_1,\theta\in\Theta_{n,r_1}}\frac{r_{n,\theta}^2\mu
_{n,\theta} k_{n,\theta,4}^{{1}/{2}}b_{n,\theta,2}^{
{1}/{2}}}{\sigma_{n,\theta}^4}<\infty \qquad\mbox{where }
b_{n,\theta,m}=E_\theta B_n^m.
\end{eqnarray}
Regarding (\ref{VnconditionedJn}) we remark that by $\mathcal{ L}_\theta(Y_{n-L_n})$ we mean
the mixture distribution $\sum_{m=n_0}^n \mathcal{ L}_\theta(Y_m)P(L_n=n-m)$,
which can be defined without requiring that $Y_{n_0},\ldots,Y_n$ and $L_n$
all be defined on the same space.
A general prescription for size biasing a sum of nonnegative variables
is given in \citet{Goldstein96}; specializing to exchangeable
indicators yields the following result.
%
%le1.1 #&#
\begin{lemma} \label{sblem}
Let $Y=\sum_{\alpha\in\mathcal{ I}}X_\alpha$ be a finite sum of
nontrivial exchangeable
Bernoulli variables $\{X_\alpha, \alpha\in\mathcal{ I}\}$, and suppose
that for $\alpha\in\mathcal{ I}$
the variables
$\{X_\beta^\alpha, \beta\in\mathcal{ I}\}$ have joint distribution
\[%\label{Xbetaalphaiscond}
\mathcal{ L}(X_\beta^\alpha, \beta\in\mathcal{ I})=\mathcal{ L}(X_\beta, \beta
\in\mathcal{ I}|X_\alpha=1).
\]
Then
\[
Y^\alpha=\sum_{\beta\in\mathcal{ I}} X_\beta^\alpha
\]
has the $Y$-size biased distribution $Y^s$, as does the mixture $Y^I$
when $I$ is a random index with values in $\mathcal{ I}$, independent of
all other variables.
\end{lemma}

\begin{pf}
First, fixing $\alpha\in\mathcal{ I}$, we show that $Y^\alpha$ satisfies
(\ref{introszb}). For given $f$,
\[
E[Y f(Y)] = \sum_{\beta\in\mathcal{ I}} E[X_\beta f(Y)]
= \sum_{\beta\in\mathcal{ I}} P[X_\beta=1]
E[f(Y)| X_\beta= 1].
\]
As exchangeability implies that $E[f(Y)| X_\beta= 1]$ does not depend
on $\beta$, we have
\[
E[Y f(Y)] =
\biggl( \sum_{\beta\in\mathcal{ I}} P[X_\beta=1] \biggr)
E[f(Y)| X_\alpha= 1]
= E[Y ]
E[f(Y^\alpha)],
\]
demonstrating the first result. The second follows easily using that
$Y^I$ is a mixture of random variables all of which have distribution
$Y^s$.
\end{pf}

Employing size bias couplings and Stein's method, \citet{ChRo10} prove
a general result to compute bounds to the normal in the Waserstein
metric. In particular, Corollary 2.2\vadjust{\goodbreak} and Construction 3A of \citet
{ChRo10} yield
%
%e14 #&#
\begin{equation}\label{ChRoWass}
d_W(\mathcal{ L}_\theta(W_{n,\theta}),\mathcal{ L}(Z)) \le0.8\frac{\mu
_{n,\theta}\Psi_{n,\theta}}{\sigma_{n,\theta}^2}+\frac{\mu
_{n,\theta
}k_{n,\theta,2}}{\sigma_{n,\theta}^3}.
\end{equation}

To compare (\ref{ChRoWass}) with one conclusion of Theorem \ref
{thmmain}, as well as to familiarize the reader with the roles of some
of the variables appearing in its formulation, we now consider its
application in the simple case where a bounded size bias coupling
exists, that is, when the bound $K_n$ on $|Y_n^s-Y_n|$ can be taken to
be a constant, say $k_n$, almost surely. In such cases we set $J_n$ to
be the empty set, and note that any constant is measurable with respect
to the trivial $\sigma$-algebra that $J_n$ generates. Conditions \hyperlink{assumptioncondlawYnvgivenK}{3} through
\hyperlink{Thetanr-compact}{6} are
easily satisfied in this case for any candidate $r_{n,\theta}$. In
particular, taking $L_n=0,s_{n,\theta}=1$ and $F_{n,\theta}=\{L_n \le
s_{n,\theta}\}$, with $J_n=\varnothing$, (\ref{VnconditionedJn}) of
Condition \hyperlink{assumptioncondlawYnvgivenK}{3} holds with $\psi
_{n,\theta
}=\theta$ and $V_n=Y_n$, and (\ref{musigKLs}) holds as $1-\mathbf{
1}(F_{n,\theta})=0$ a.s. As $(n-L_n,\psi_{n,\theta})=(n,\theta)$,
Condition \hyperlink{unifvarbnd}{4} holds with $c_1=c_2=1$. As $V_n=Y_n$ we
may take $B_n=0$ in Condition \hyperlink{assumptionYnmYnVnprime}{5}, and as
$L_n=0$ Condition \hyperlink{Thetanr-compact}{6}a is satisfied. Hence, only
Conditions \hyperlink{Psinisrootn}{1} and \hyperlink{assumptionKn}{2} are in force, and
Theorem~\ref{thmmain} obtains with
\[
r_{n,\theta}^{-1}=\frac{\mu_{n,\theta} \Psi_{n,\theta}}{\sigma
_{n,\theta
}^2} + \frac{\mu_{n,\theta}k_n^2}{\sigma_{n,\theta}^3},
\]
yielding a Kolmogorov bound that, up to constants, agrees with the
Wasserstein bound (\ref{ChRoWass}) in this particular case.

Bounded size bias couplings exist when $Y_n$ is the sum of independent,
bounded nonnegative random variables, or a sum of bounded, nonnegative
locally dependent variables with bounded dependence neighborhood sizes,
as studied, for instance, in \citet{Go05}. In addition, bounded size
bias couplings can also be constructed in cases of global dependence;
see \citet{GoZh} or \citet{Goldstein10}.

We next apply Theorem~\ref{thmmain} to vertex degree counts in the
Erd\H{o}s--R\'enyi random graph. The proof of Theorem~\ref{thmmain} is
given in Section~\ref{secproofthmmain}.

%s2 #&#
\section{Vertex degree in the Erd\H{o}s--R\'enyi random graph}
We apply Theorem~\ref{thmmain} to bound the error in the normal
approximation to the distribution of the number of vertices of a given
degree in the Erd\H{o}s--R\'enyi \citeyearpar{Erdos59} random graph
$G_n$; see also \citet{Bollobas85}. With $n \in\mathbb{N}$ we take the
vertex set of $G_n$ to be $\mathcal{ I}_n=\{1,\ldots,n\}$, and the
indicators $\xi_{u,v}$ of the presence of edges between distinct
vertices $u$ and $v$ to be independent Bernoulli variables with a
common success probability. No vertex is connected to itself, and we
set $\xi_{u,u}=0$ for all $u \in\mathcal{ I}_n$.

The number $Y_n$ of vertices of degree $d$ of $G_n$ has been the object
of much study. For a sequence of graphs with connectivity probability
$p$ depending on $n \in\mathbb{N}$,
\citet{Karonski87}
proved the asymptotic normality of $Y_n$ when $n^{(d+1)/d} p\rightarrow
\infty$ and $n p\rightarrow0$, or $n p \rightarrow\infty$ and $n
p-\log n
-d \log\log n \rightarrow- \infty$; see also \citet{Palka84} and
\citet{Bollobas85}.
Asymptotic normality of $Y_n$ when $np \rightarrow c > 0$ was obtained
by \citet{Barbour89}, and \citet{Kordecki90}
for nonsmooth functions of $Y_n$ in the case $d=0$. \citet
{Neammanee09} obtain a Kolmogorov distance bound between $Y_n$ and the
normal with rate $n^{-1/2+\varepsilon}$ for all $\varepsilon>0$ when $\operatorname
{Var}(Y_n)$ is of order $n$. Other univariate
results on asymptotic normality of counts on random graphs
are given in \citet{Janson91}, and references therein. \citet{Goldstein96}
obtain smooth function bounds for the vector whose $k$ components count the
number of vertices of fixed degrees $d_1,d_2,\ldots,d_k$ when
$p=\theta
/(n-1) \in(0,1)$ for fixed $\theta$, implying
asymptotic multivariate joint normality.

We focus on the counts of vertices of some fixed degree $d \in\mathbb
{N}$, the case $d=0$ of isolated vertices having already been handled
by \citet{Kordecki90}. Set
%
%e15 #&#
\begin{equation}\label{defThetangraph}
\Theta_n=(0,n-1) \cap(0,b]\qquad \mbox{for all $n \ge d+1$}
\end{equation}
with $b$ some arbitrarily large constant, and let the connectivity
probability between the vertices of $G_n$ be given by $\theta/(n-1)$
for $n \ge d+1, \theta\in\Theta_n$. For $v \in\mathcal{ I}_n$ let
\[%\label{defdeg}
D_n(v)=\sum_{w \in\mathcal{ I}_n}\xi_{v,w},\qquad   X_{n,v}=\mathbf{
1}\bigl(D_n(v)=d\bigr) \quad\mbox{and}\quad
Y_n=\sum_{v \in\mathcal{ I}_n}X_{n,v},
\]
the degree of vertex $v$, the indicator that vertex $v$ has degree $d$,
and the number of vertices of degree $d$ of $G_n$, respectively.

From \citet{Goldstein96}, for all $n \ge d+1$ and $\theta\in\Theta
_n$, the mean $\mu_{n,\theta}$ and variance $\sigma_{n,\theta}^2$ of
$Y_n$ are given explicitly by
%
%e16 #&#
\begin{equation}
\label{cov}
\mu_{n,\theta}=n\tau_{n,\theta}\quad  \mbox{and}\quad  \sigma
_{n,\theta
}^2= n\tau_{n,\theta}^2
\biggl[\frac{(d-\theta)^2}{\theta(1-\theta/(n-1))}-1\biggr]
+ n\tau_{n,\theta},
\end{equation}
where
%
%e17 #&#
\begin{equation}
\label{deftauntheta}
\tau_{n,\theta}=\pmatrix{{n-1} \vspace*{2pt}\cr d}\biggl( \frac{\theta}{n-1}\biggr)^d
\biggl(1-\frac{\theta}{n-1}\biggr)^{n-1-d}.
\end{equation}

%th2.1 #&#
\begin{theorem}
\label{thmdegree} For any $d \in\mathbb{N}$ and $b>0$ there exists a
constant $C$ such that for all $n \ge d+1$ and all $\theta\in\Theta
_n$ given in (\ref{defThetangraph}), the normalized count
$W_{n,\theta}$ in (\ref{defWntheta}) of the number $Y_n$ of vertices
with degree $d$ in the Erd\H{o}s--R\'enyi random graph $G_n$ on~$n$
vertices, with edges connecting each distinct pair independently with
probability $\theta/(n-1)$, satisfies
\[
\sup_{z \in\mathbb{R}}|P_\theta(W_{n,\theta} \le z)-P(Z \le z)|
\le
C/r_{n,\theta}  \qquad\mbox{for all $n \ge d+1$},
\]
where $Z$ is a standard normal variable and
%
%e18 #&#
\begin{equation}\label{defrtau}
r_{n,\theta}=\sqrt{n \tau_\theta}\qquad \mbox{with } \tau_\theta
=e^{-\theta
}\theta^d/d!.
\end{equation}
\end{theorem}

By applying Stein's method, \citet{Barbour89} obtain a bound of order
$1/\sqrt{n \tau_{n,\theta}}$ in the metric
$d_L$ defined as the supremum over Lipschitz functions
%
%e19 #&#
\begin{equation}\label{defd1}
d_L(\mathcal{ L}(X),\mathcal{ L}(Y)) = \sup_h \frac{|Eh(X)-Eh(Y)|}{\Vert h\Vert +\Vert h'\Vert }.
\end{equation}
As Lemma~\ref{degreetaunconverges} shows that $\tau_{n,\theta
}/\tau
_\theta$ converges uniformly to 1 over $\Theta_n$, the Kolmogorov bound
of order $1/\sqrt{n \tau_\theta}$ provided by Theorem~\ref{thmdegree}
is of the same order as the $d_L$ bound. As remarked in \citet
{Barbour89}, a bound of size $\varepsilon_n$ in the $d_L$ metric yields a
bound in the Kolmogorov metric of order $O(\varepsilon_n^{1/2})$, which
can at times be improved to $O(\varepsilon_n)$ ``at the cost of much
greater effort.''

Though we do not cover the case $d=0$ of isolated vertices, handled in
\citet{Kordecki90}, our proof can be extended to apply there by
appending additional arguments that are separate, but similar to, those
for the case $d \in\mathbb{N}$. Note, for example, the difference in
the behavior of the function $\tau_\theta$ at zero for these two ranges
of $d$.

Following Lemma~\ref{sblem} for the case of vertex degrees yields a
coupling where for each $n \ge d+1$ and vertex $v \in\mathcal{ I}_n$ one
constructs a graph $G_n^v$ from $G_n$ having the distribution of $G_n$
conditioned on $X_{n,v}=1$, or equivalently, on $D_n(v)=d$; this
coupling has previously been applied by \citet{Barbour92} and \citet
{Goldstein96}. The graph $G_n^v$ is obtained from $G_n$ by adding or
removing edges of $v$ as needed. Mixing over $v$ as indicated by Lemma
\ref{sblem} yields a variable $Y_n^s$ having the $Y_n$-size bias distribution.

In the course of constructing $G_n^v$ one also obtains a set $\mathcal{
R}_n^v$ holding the collection of vertices other than $v$ that are
affected by the size bias operation. In particular, if $D_n(v)=d$, then
$G_n^v=G_n$ and $\mathcal{ R}_n^v=\varnothing$. If $D_n(v)>d$, then $G_n^v$
is formed by removing from $G_n$ the edges between $v$ and the vertices
in the subset $\mathcal{ R}_n^v$ of neighbors $\{u\dvtx  \xi_{u,v}=1\}$ of $v$,
chosen with uniform conditional distribution given $G_n$ over all
subsets of the neighbors of $v$ of size $D_n(v)-d$. Similarly, if
$D_n(v)<d$, then $G_n^v$ is formed by adding edges to $G_n$ between $v$
and vertices in $\mathcal{ R}_n^v$, chosen with uniform conditional
distribution given $G_n$ over all subsets of the nonneighbors $\{u\dvtx  u
\not= v, \xi_{u,v}=0\}$ of $v$ of size $d-D_n(v)$.

Now let $X_{n,w}^v$ be the indicator that vertex $w$ has degree $d$ in
$G_n^v$ and
\[
Y_n^v=\sum_{w \in\mathcal{ I}_n} X_{n,w}^v,
\]
the number of degree $d$ vertices in $G_n^v$. When $I_n$ is chosen
uniformly over $\mathcal{ I}_n$, independent of all other variables, Lemma
\ref{sblem} yields that $Y_n^s=Y_n^{I_n}$ has the $Y_n$-size biased
distribution.
Similarly setting $\mathcal{ R}_n^s=\mathcal{ R}_n^{I_n}$, all vertices not in
$\{I_n\} \cup\mathcal{ R}_n^s$ have the same degree in both $G_n$ and
$G_n^s$, and as $I_n \notin\mathcal{ R}_n^s$, letting
%
%e20 #&#
\begin{equation}\label{InsizeRns}
\mathcal{ A}_n= \{I_n\} \cup\mathcal{ R}_n^s\qquad \mbox{we have } |\mathcal{ A}_n|= 1+
|d-D_n(I_n)|.
\end{equation}

We prove Theorem~\ref{thmdegree} by verifying the hypotheses of
Theorem~\ref{thmmain} for the
size bias construction just given. With $\tau_{n,\theta}$ as in (\ref
{deftauntheta}), and recalling (\ref{cov}), let
%
%e21 #&#
\begin{equation}\label{degreedefdelta}\qquad
\delta_{n,\theta}=\tau_{n,\theta}\biggl[\frac{(d-\theta
)^2}{\theta
(1-\theta/(n-1))}-1\biggr]+1
\qquad\mbox{so that }
\sigma_{n,\theta}^2 = n \tau_{n,\theta} \delta_{n,\theta},
\end{equation}
and correspondingly, with $\tau_\theta$ as in (\ref{defrtau}), let
%
%e22 #&#
\begin{equation}\label{defdelta}
\delta_\theta=
\tau_\theta\biggl[\frac{(d-\theta)^2}{\theta}-1\biggr]+1.
\end{equation}
%
%We say that a sequence of functions $f_{n,\theta}$ converges uniformly
%to $f_\theta$ on $\Theta_n$ if
%f_{n,\theta}-f_\theta| | =0.

With the help of a technical lemma placed at the end of this section,
we present the proof of Theorem~\ref{thmdegree}. Throughout we let
$C_j$ denote a constant not depending on $n$ or $\theta$, and not
necessarily the same at each occurrence.

\begin{pf*}{Proof of Theorem \protect\ref{thmdegree}}
Let $n_0=d+1$. For $n \ge n_0$ and $\theta\in\Theta_n$ the
binomial and Poisson probabilities $\tau_{n,\theta}$ and $\tau
_\theta$
in (\ref{deftauntheta}) and (\ref{defrtau}), respectively, lie in
$(0,1)$, and hence $\sigma_{n,\theta}^2$ of (\ref{cov}) and
$r_{n,\theta
}$ are positive for all such $n$ and $\theta$. Let $r_1>0$ be
arbitrary. In place of naming $n_1$ explicitly, we show the remaining
conditions of Theorem~\ref{thmmain} are satisfied for all $n$
sufficiently large. Since $r_{n,\theta} \le\sqrt{n}$ inequality~(\ref{bdawayzerornt}) holds for any $n_1 \ge n_0$.

From \citeauthor{Chen10} [(\citeyear{Chen10}), equation (12.17)], following \citet{Goldstein96}, for $Y_n^s$
having the $Y_n$-size biased distribution as constructed above, we obtain
\[%\label{graphPsibd}
\Psi_{n,\theta}^2 \le
\C n^{-1}( 24 \theta+ 48\theta^2 + 144\theta^3
+48d^2+144\theta
d^2+12)
\]
and hence
\[
\sup_{\theta\in\Theta_n} \Psi_{n,\theta}\le\frac{\C}{\sqrt{n}}.
\]

To complete the verification of Condition \hyperlink{Psinisrootn}{1}, Lemma
\ref
{degreetaunconverges} gives that over $\Theta_n$ the ratio $\delta
_\theta/\delta_{n,\theta}=\delta_\theta\mu_{n,\theta}/\sigma
_{n,\theta
}^2$ converges uniformly to $1$, and $\delta_\theta$ in (\ref
{defdelta}) is bounded away from zero. Hence for all $n$ sufficiently
large and all $\theta\in\Theta_n$, we have
%
%e23 #&#
\begin{equation}\label{degreemuoversigma2}
\frac{\mu_{n,\theta}}{\sigma_{n,\theta}^2} \le\frac{2}{\delta
_\theta}
\le\C\quad
\mbox{and so}\quad \frac{r_{n,\theta} \mu_{n,\theta} \Psi_{n,\theta
}}{\sigma
_{n,\theta}^2} \le
\Cl{Psi.tau.theta}\sqrt{\tau_\theta} \le\Cr{Psi.tau.theta}
\end{equation}
as $\tau_\theta\le1$ for all $\theta\in\Theta_n$.

Turning to Condition \hyperlink{assumptionKn}{2}, let
\[%\label{defAnJnFn}
J_n=(I_n, D_n(I_n))\quad \mbox{and}\quad\mathcal{ F}_n=\sigma\{J_n \};
\]
that is, $\mathcal{ F}_n$ is the $\sigma$-algebra generated by the chosen
vertex and its degree. Further, let
\[%\label{degreedefKntheta}
K_n=1+d+D_n(I_n).
\]
Clearly $K_n$ is $\mathcal{ F}_n$-measurable, and recalling that vertices
not in $\mathcal{ A}_n$ of (\ref{InsizeRns}) have the same degree in
both\vadjust{\goodbreak}
$G_n$ and $G_n^s$, taking the difference between $Y_n^s$ and $Y_n$ yields
\[%\label{defYnsYdeg}
Y_n^s-Y_n=\sum_{w \in\mathcal{ A}_n}(X_{n,w}^{I_n}-X_{n,w}),
\]
and (\ref {InsizeRns}) yields
\[
|Y_n^s-Y_n| = 1+ |d-D(I_n)| \le K_n.
\]

Next, for all $m \in\mathbb{N}$ we have
%
%e24 #&#
\begin{equation}\label{Km2mbound}
K_n^m \le2^{m-1}\bigl( (1+d)^m +D_n(I_n)^m \bigr).
\end{equation}
To bound the moments of $K_n$, using \citet{Riordan37} for the first
equality below, with $S_{j,m}$ the Stirling numbers of the second kind,
$(n)_j$ the falling factorial, $\Cl{sumdeg}{_{,m}} = m\max_{1 \le j
\le
m}S_{j,m}$ and $D \sim\operatorname{Bin}(n-1,p)$, we obtain
\begin{eqnarray*}%\label{binmom}
E D^m &=& \sum_{j=1}^m S_{j,m} (n-1)_j p^j \le\sum_{j=1}^m S_{j,m}
(n-1)^j p^j\\
&\le&\Cr{sumdeg}{_{,m}} \bigl((n-1)p+(n-1)^mp^m\bigr).
\end{eqnarray*}
In particular $E_\theta D_n(v)^m \le\Cr{sumdeg}{_{,m}} (b+b^m)$, and
as $D_n(I_n)$ is the mixture of the identical distributions $D_n(v)$
over $v \in\mathcal{ I}_n$, it obeys the same upper bound.
Taking expectation in (\ref{Km2mbound}), we find that there exists
constants $\Cl{K}{_{,m}},m \in\mathbb{N}$ such that
%
%e25 #&#
\begin{eqnarray}\label{degreeknthetambounded}
k_{n,\theta,m} \le\Cr{K}{_{,m}} \qquad\mbox{for all $n \in\mathbb
{N},\theta \in\Theta_n$ and $m \in\mathbb{N}$.}
\end{eqnarray}
Now, using (\ref{degreeknthetambounded}) for the first inequality in
(\ref{degreemusigma3}), the first inequality in (\ref
{degreemuoversigma2}) for the second inequality, the second equality
of (\ref{degreedefdelta}) for the first equality, and Lemma \ref
{degreetaunconverges} both to obtain the third inequality, and the
boundedness of $\delta_\theta$ away from zero for the fourth, we obtain
that for all $n$ sufficiently large and $\theta\in\Theta_n$,
%
%e26 #&#
\begin{equation}\label{degreemusigma3}
\qquad \frac{r_{n,\theta}\mu_{n,\theta}k_{n,\theta,4}^{1/2}}{\sigma
_{n,\theta
}^3} \le\frac{\Cr{K}{_{,4}}^{1/2}r_{n,\theta}\mu_{n,\theta
}}{\sigma
_{n,\theta}^3}
\le\frac{\Cl{C: (20)} r_{n,\theta}}{\sigma_{n,\theta}} =
 \frac
{\Cr{C: (20)}
\sqrt{\tau_\theta}}{\sqrt{\tau_{n,\theta}\delta_{n,\theta}}}
\le\frac
{\C}{\sqrt{\delta_\theta}} \le\C.
\end{equation}
Hence inequality (\ref{kn-sum2-conv}), sufficient for (\ref
{assumptionKninequality}), is satisfied, and Condition \hyperlink{assumptionKn}{2} holds.

Turning to Condition \hyperlink{assumptioncondlawYnvgivenK}{3}, for $n \ge
d+2$, let
%
%e27 #&#
\begin{equation}\label{defpsitheta}
 L_n=1, \qquad  s_{n,\theta}=1, \qquad \psi_{n,\theta}=
\biggl( \frac{n-2}{n-1}\biggr)\theta\quad\mbox{and}\quad F_{n,\theta}=\{L_n
\le s_{n,\theta}\},\hspace*{-35pt}
\end{equation}
and note therefore that conditions holding on $F_{n,\theta}$ must hold
on the entire probability space. Clearly $L_n$ takes values in $\{
0,1,\ldots,n\}$ as required and $n-s_{n,\theta} \ge n_0$ for any $n
\ge
d+2$. Being constants, $L_n$ and $\psi_{n,\theta}$ are $\mathcal{ F}_n$
measurable, hence $F_{n,\theta} \in\mathcal{ F}_n$. By (\ref
{defpsitheta}) and $\theta\in\Theta_n$ we have that $\psi
_{n,\theta}
\in(0,b] \cap(0,n-2) = \Theta_{n-1}=\Theta_{n-L_n}$, verifying the
first part of (\ref{VnconditionedJn}).

Regarding the second part of (\ref{VnconditionedJn}), let $H_n$ be
the graph $G_n$ with the vertex $I_n$ and its incident edges removed,
relabeling the remaining vertices $\{1,\ldots,n-1\}$ by preserving
their relative order. Let $V_n$ be the number of degree $d$ vertices of
$H_n$. By counting the number of degree $d$ vertices, the
distributional equality in (\ref{VnconditionedJn}) is a consequence of
%
%e28 #&#
\begin{equation}\label{GnInCondisGnminusK}
\mathcal{ L}_\theta(H_n|I_n, D_n(I_n)) = \mathcal{ L}_{\psi_{n,\theta}}(G_{n-1}).
\end{equation}
The graph $H_n$ is determined by $\{\xi_{u,v}\dvtx  \{u,v\} \subset\mathcal{
I}_n \setminus\{I_n\}\}$, which is independent of the $\sigma$-algebra
generated by $\{I_n, \xi_{I_n,v}, v \in\mathcal{ I}_n\}$, with respect to
which $I_n$ and $D_n(I_n)$ are measurable. Hence $H_n$ is independent
of the conditioning event in (\ref{GnInCondisGnminusK}), and therefore
its conditional and unconditional distribution agree. In particular,
conditional on $\{I_n,D_n(I_n)\}$, the edge indicators of $H_n$ are
independent with common success probability
\[
\frac{\theta}{n-1}= \frac{\psi_{n,\theta}}{n-2},
\]
so (\ref{GnInCondisGnminusK}) holds. Inequality (\ref{musigKLs})
holds trivially, as $P(L_n > 1)=0$. Hence Condition~\hyperlink{assumptioncondlawYnvgivenK}{3} holds.

By Lemma~\ref{degreetaunconverges}, Condition \hyperlink{unifvarbnd}{4}
holds with $c_1=c_2=2$.

Regarding Condition \hyperlink{assumptionYnmYnVnprime}{5}, as only the degrees
of vertex $I_n$ and its neighbors are different in the graphs $G_n$ and
$H_n$, we have
\[
|Y_n-V_n| \le1+ D(I_n) \le K_n,
\]
and we set $B_n=K_n$, so $\mathcal{ F}_n$-measurable. We now finish the
verification of Condition~\hyperlink{assumptionYnmYnVnprime}{5} by showing
(\ref
{inequalityYnmYnVnprimesufficient}), sufficient for (\ref
{inequalityYnmYnVnprime}), is satisfied. By (\ref
{degreeknthetambounded}), that $\mu_{n,\theta}=n\tau_{n,\theta}$ and
the second equality in (\ref{degreedefdelta}), for all $n$
sufficiently large and all $\theta\in\Theta_n$, we have
\[
\frac{r_{n,\theta}^2\mu_{n,\theta} k_{n,\theta,4}^{
{1}/{2}}b_{n,\theta,2}^{{1}/{2}}}{\sigma_{n,\theta}^4} \le\frac
{\tau
_\theta(\Cr{K}{_{,4}}\Cr{K}{_{,2}})^{1/2}}{\tau_{n,\theta} \delta
_{n,\theta}^2} \le\C,
\]
where the final inequality follows from Lemma \ref
{degreetaunconverges}, yielding that $\tau_{n,\theta}/\tau_\theta$
and $\delta_{n,\theta}/\delta_\theta$ converge uniformly to 1 on
$\Theta
_n$, and that $\delta_\theta$ is bounded away from zero on $(0,b]$.

Lastly, Condition \hyperlink{Thetanr-compact}{6}a holds with $l_{n,0}=1$ for
all $n \ge d+2$, completing the verification of all conditions of
Theorem~\ref{thmmain}.
\end{pf*}

The proof of Lemma~\ref{degreetaunconverges} is straightforward, and
is therefore omitted.

%For any $b>0$ the function
%e^\theta(1-\frac{\theta}{n} )^n
%converges uniformly to 1 on the interval $(0,b]$.

%(-\infty,0]$, and that for all $x \in(0,1)$ a %first order Taylor
%expansion about zero yields
%Hence, for all $\theta\in(0,b]$,
%| (1-\frac{\theta}{n})^n-e^{-\theta} | \le n
%| \log(1-\frac{\theta}{n}) + %\frac{\theta}{n} |
%so the left-hand side tends to zero uniformly as $n \rightarrow
%| e^\theta(1-\frac{\theta}{n} )^n - 1| = e^

%Recall the definitions of $\mu_{n,\theta}, \sigma_{n,\theta}^2$ and $

%le2.1 #&#
\begin{lemma}
\label{degreetaunconverges}
With $\tau_{n,\theta}, \tau_\theta, \delta_{n,\theta}$ and
$\delta
_\theta$ given by (\ref{deftauntheta}), (\ref{defrtau}), (\ref
{degreedefdelta}) and (\ref{defdelta}), respectively, for all $d
\in
\mathbb{N}$ and all $b>0$ the function $\delta_\theta$
is bounded away from zero and infinity over $(0,b]$, and the ratios
\begin{eqnarray*}%\label{degreelistratios}
\frac{\tau_{n,\theta}}{\tau_\theta},\qquad
\frac{\delta_{n,\theta}}{\delta_\theta},\qquad
\frac{r_{n,\theta}}{r_{n-1,\psi_{n,\theta}}}\quad
\mbox{and}\quad
\frac{\sigma_{n,\theta}^2}{\sigma_{n-1,\psi_{n,\theta}}^2}
\end{eqnarray*}
and their reciprocals converge uniformly to 1 on $(0,b]$ as $n$ tends
to infinity.\vadjust{\goodbreak}
\end{lemma}

\section{\texorpdfstring{Proof of Theorem \protect\ref{thmmain}}{Proof of Theorem 1.1}}
\label{secproofthmmain}
\resetconstant
We begin the proof of Theorem~\ref{thmmain} with the following lemma.
%
%le3.1 #&#
\begin{lemma}
\label{big-recursion}
Suppose that for some $n_1 \in\mathbb{N}_0$ the nonnegative numbers~$f$,
$\{p_{n,l}\}_{n\ge n_1,0 \le l \le n}$ and $\{a_n\}_{n \ge0}$ satisfy
%
%e29 #&#
\begin{eqnarray}\label{anpnrec}
&&a_n \le\sum_{l=0}^n a_{n-l}p_{n,l} + f \qquad\mbox{for all $n \ge n_1$\quad
and}
\nonumber
\\[-8pt]
\\[-8pt]
\nonumber
&&\qquad\tau\in(0,1)\qquad \mbox{where } \tau=\sup_{n\ge n_1}\sum_{l=0}^n p_{n,l}.
\end{eqnarray}
Then $\sup_{n \ge0} a_n < \infty$.
\end{lemma}

\begin{pf} As for all $n \ge n_1$ we have $p_{n,0} \le\tau< 1$, letting
\[
q_{n,l}=\frac{p_{n,l}}{1-p_{n,0}}\qquad \mbox{for $1 \le l \le n$\quad and}\quad a =
\frac{f}{1-\tau},
\]
(\ref{anpnrec}) implies
\begin{eqnarray*}
a_n \le\sum_{l=1}^n a_{n-l}q_{n,l} + a\qquad \mbox{with }
0 \le\sum_{l=1}^n q_{n,l} \le\frac{\tau-p_{n,0}}{1-p_{n,0}} \le
\tau
\mbox{ for all $n \ge n_1$.}
\end{eqnarray*}

Letting $\alpha=\max_{0 \le n \le n_1}a_n$ and $c=\max\{a,\alpha
(1-\tau
)\}$, the sequence $\{b_n\}_{n \ge0}$ defined by
\[
b_n=\alpha\qquad \mbox{for $0 \le n \le n_1$\quad and}\quad   b_{n+1}=\tau
b_n + c\qquad  \mbox{for $n \ge n_1$}
\]
has, for $n \ge n_1$, the explicit form
\[
b_n=\gamma\tau^{n-n_1} + \frac{c}{1-\tau}\qquad \mbox{where }
\gamma=\alpha-\frac{c}{1-\tau}.
\]
Since $\gamma\le0$ and $\tau\in(0,1)$, the sequence $\{b_n\}_{n \ge
0}$ is nondecreasing with limit $c/(1-\tau)$, and hence is bounded. We
complete the proof by showing that for all $n \in\mathbb{N}_0$ we have
$a_m \le b_m$ for all $0 \le m \le n$. Clearly the statement holds for
$0 \le n \le n_1$. Assuming it true for some $n \ge n_1$, using the
induction hypotheses, the definition of $c$ and that $b_n$ is nondecreasing,
\begin{eqnarray*}
a_{n+1} &\le&\sum_{l=1}^{n+1} a_{n+1-l}q_{n+1,l} + a \le\sum
_{l=1}^{n+1} b_{n+1-l}q_{n+1,l} + c \le b_n\sum_{l=1}^{n+1} q_{n+1,l} +
c\\
 &\le&\tau b_n + c = b_{n+1}.
\end{eqnarray*}
\upqed\end{pf}

The following proof is based on the inductive argument of
\citet{Bolthausen84}.

\begin{pf*}{Proof of Theorem \protect\ref{thmmain}} With $r \ge0$, recall that
$\Theta_{n,r}=\{\theta\in\Theta_n\dvtx \break r_{n,\theta} \ge r\}$, and let
%
%e30 #&#
\begin{eqnarray}\label{defdeltanr}
\delta(n,r)=
\sup_{z \in\mathbb{R}, \theta\in\Theta_{n,r}} |P_\theta
(W_{n,\theta}
\le z)-P(Z \le z)| \qquad\mbox{for $n \ge n_0$.}
\end{eqnarray}
First note that (\ref{thmmainstatement}) of Theorem~\ref{thmmain}
can be made to hold
whenever $r_{n,\theta} < r_1$ by taking $C \ge r_1$. By (\ref
{bdawayzerornt}) the cases $n_0 \le n < n_1$ and $r_{n,\theta} \ge
r_1$ can be handled in this same manner. Hence it suffices to show that
there exists some $C$ such that
%
%e31 #&#
\begin{eqnarray}\label{pfgoal}
\delta(n,r) \le C/r\qquad \mbox{for $n \ge n_1$ and $r \ge r_1$.}
\end{eqnarray}

For $z \in\mathbb{R}$ and $\lambda>0$ let $h_{z,\lambda}$ be the
smoothed indicator
\[
h_{z,\lambda}(x)=\cases{
1, & \quad$x \le z,$ \vspace*{2pt}\cr
1+(z-x)/\lambda,& \quad$z< x \le z+\lambda,$\vspace*{2pt}\cr
0, & \quad $x > z+\lambda$}
\]
and let $Nh_{z,\lambda}=Eh_{z,\lambda}(Z)$ with $Z$ a standard normal
variable.
Let $f(x)$ be the unique bounded solution to the Stein equation for
$h_{z,\lambda}(x)$ [see, e.g., \citet{Chen10}]
%
%e32 #&#
\begin{equation}\label{eqSteinf}
h_{z,\lambda}(x)-Nh_{z,\lambda} = f'(x)-xf(x).
\end{equation}

Let $n \ge n_1, \theta\in\Theta_{n,r}$ for some $r \ge r_1$, $z \in
\mathbb{R}$ and $\lambda>0$. Recalling
$W_{n,\theta}=(Y_n-\mu_{n,\theta})/\sigma_{n,\theta}$, with a slight
abuse of notation, set
\[
W_{n,\theta}^s=\frac{Y_n^s-\mu_{n,\theta}}{\sigma_{n,\theta}}.
\]

Substituting $W_{n,\theta}$ for $x$ in (\ref{eqSteinf}) and taking
expectation, and dropping the subscript $\theta$ when not essential
below, we obtain
%
%e33 #&#
\begin{equation}
\label{h-diff}
E_\theta h_{z,\lambda}(W_n)-Nh_{z,\lambda} = E_\theta[f'(W_n)-W_n f(W_n)].
\end{equation}
Beginning with the second term on the right-hand side of (\ref
{h-diff}), from the definition of $W_{n,\theta}$ and the size bias
relation (\ref{introszb}),
we have
\begin{eqnarray*}
E_\theta[W_nf(W_n)]=\frac{1}{\sigma_n}E_\theta[(Y_n-\mu
_n)f(W_n)]=\frac{\mu_n}{\sigma_n}E_\theta\bigl(f(W_n^s)-f(W_n)
\bigr).
\end{eqnarray*}
Taking absolute value and applying the triangle inequality, we obtain
%
%e34 #&#
\begin{eqnarray}\label{two-main-terms}
&&|E_\theta h_{z,\lambda}(W_n)-Nh_{z,\lambda}|\nonumber\\
&&\qquad= |E_\theta[f'(W_n)-W_nf(W_n)] |\nonumber\\
&&\qquad= \biggl| E_\theta\biggl[f'(W_n)-\frac{\mu_n}{\sigma_n}
\bigl(f(W_n^s)-f(W_n) \bigr)\biggr]\biggr|\nonumber\\
&&\qquad= \frac{\mu_n}{\sigma_n} \biggl|E_\theta\biggl[\frac
{\sigma
_n}{\mu_n}f'(W_n) -\bigl (f(W_n^s)-f(W_n) \bigr)\biggr]\biggr|
\nonumber
\\[-8pt]
\\[-8pt]
\nonumber
 &&\qquad= \frac{\mu_n}{\sigma_n} \biggl| E_\theta
\biggl[ \biggl( \frac{\sigma_n}{\mu_n} -(W_n^s-W_n)\biggr) f'(W_n)
+(W_n^s-W_n)f'(W_n)
\\
&&\hspace*{175pt}\qquad{}- \bigl(f(W_n^s)-f(W_n) \bigr)\biggr]\biggr |
\nonumber\\
&&\qquad\le \frac{\mu_n}{\sigma_n} \biggl| E_\theta\biggl[
\biggl(
\frac
{\sigma_n}{\mu_n} -(W_n^s-W_n)\biggr) f'(W_n)\biggr] \biggr|
\nonumber\\
&&\qquad\quad{}+\frac{\mu
_n}{\sigma_n} \biggl| E_\theta\biggl[ \int_0^{W_n^s-W_n}[f'(W_n)
-f'(W_n+t)]\,dt \biggr]\biggr|.\nonumber
\end{eqnarray}

From the size bias relation (\ref{introszb}) with $f(x)=x$, we obtain
$\mu_n E_\theta[Y_n^s]=E_\theta[Y_n^2]$, and therefore
%
%e35 #&#
\begin{equation}
\label{ExpSnV-Sn}\qquad
E_\theta[W_n^s-W_n]=E_\theta\biggl[\frac{Y_n^s-Y_n}{\sigma_n}
\biggr]=\frac{1}{\sigma_n}\biggl[ \frac{E_\theta Y_n^2}{\mu_n}- \mu_n
\biggr]=\frac{1}{\sigma_n \mu_n}\sigma_n^2=\frac{\sigma_n}{\mu_n}.
\end{equation}
Now applying (\ref{ExpSnV-Sn}) and
$|f'(x)| \le1$ from \citeauthor{Chen04} [(\citeyear{Chen04}), equation~(4.6)] [see also \citet
{Chen10}, Lemma~2.5], by conditioning on $W_n$ the first term of (\ref
{two-main-terms}) may be bounded by
%
%e36 #&#
\begin{eqnarray}
\label{term1}
&&\frac{\mu_n}{\sigma_n}\biggl | E_\theta\biggl[ E_\theta
\biggl(\frac
{\sigma_n}{\mu_n} -(W_n^s-W_n)\Big\vert W_n\biggr) f'(W_n)
\biggr]\biggr|
\nonumber
\\[-8pt]
\\[-8pt]
\nonumber
&&\qquad
\le\frac{\mu_n}{\sigma_n}\sqrt{\operatorname{Var}E_\theta
(W_n^s-W_n|W_n)} =
\frac{\mu_n}{\sigma_n^2}\Psi_n,
\end{eqnarray}
recalling the definition of $\Psi_n$ in (\ref{defpsin}).

Moving now to the second term of (\ref{two-main-terms}),
\citeauthor{Bolthausen84} [(\citeyear{Bolthausen84}), equation~(2.4)] gives
\[
|f(x)|\le1\quad\mbox{and}\quad|xf(x)|\le1,
\]
and combining these inequalities with $|f'(x)| \le1$ and (\ref
{eqSteinf}) as in \citeauthor{Bolthausen84} [(\citeyear{Bolthausen84}), equation~(2.5)] yields
\begin{eqnarray*}%\label{ineq-bolt}
|f'(x)-f'(x+t)|\le|t|\biggl(1+|x|+\frac{1}{\lambda}\int_0^1
1_{[z,z+\lambda]}(x+ut)\,du \biggr).
\end{eqnarray*}
Hence, applying the bound $|Y_n^s-Y_n| \le K_n$, the second term in
(\ref{two-main-terms}) may be bounded by
%
%e37 #&#
\begin{equation}\label{from-bolt}
\frac{\mu_n}{\sigma_n}E_\theta\int_{-K_n/\sigma_n}^{K_n/\sigma
_n} |t|
\biggl( 1+|W_n|+\frac{1}{\lambda}\int_0^1 1_{[z,z+\lambda]}(W_n+ut)\,du
\biggr)\,dt,
\end{equation}
yielding three terms.

For the first two terms in (\ref{from-bolt}) we obtain
%
%e38 #&#
\begin{equation}\label{term2}
\frac{2\mu_n}{\sigma_n} E_\theta\biggl( (1+|W_n|) \int
_0^{K_n/\sigma
_n} t \,dt\biggr) = \frac{\mu_n}{\sigma_n^3}E_\theta[(1+|W_n|) K_n^2].
%&\le& \frac{\mu_n}{\sigma_n}(E_\theta[W_n^s-W_n]^4)^{1/2} \le\frac{
\end{equation}
%
%%% can probably eliminate fourth moment by conditioning on W_n,
%pulling out Y_n^s-Y_n, express in terms of \psi_n
%%% leaving E_\theta[|S_n| | \mathcal{ F}_n] \le\sqrt{E_\theta[S_n^2 | {
%%% may have to go back to Bolthausen's bound to see about removing
%absolute value

Next, as $|t| \le K_n/\sigma_n$ in the region of integration, we may
bound the expectation of the remaining term in (\ref{from-bolt}) by
%
%e39 #&#
\begin{equation}\label{befdec}
%&\le&
\frac{\mu_n}{\lambda\sigma_n^2} E_\theta\biggl( K_n \int
_{-K_n/\sigma
_n}^{K_n/\sigma_n} \int_0^1 \mathbf{1}_{[z,z+\lambda]}(W_n+ut)\,du\, dt
\biggr).
\end{equation}
Clearly,
%
%e40 #&#
\begin{equation}
\label{indicator2terms}
\mathbf{1}_{[z,z+\lambda]}(W_n+ut)
\le(1-\mathbf{1}_{F_{n,\theta}})+\mathbf{1}_{[z,z+\lambda]}(W_n+ut)\mathbf{
1}_{F_{n,\theta}}.
\end{equation}
Substituting (\ref{indicator2terms}) into (\ref{befdec}), the first
term in (\ref{indicator2terms}) gives rise to the expression
%
%e41 #&#
\begin{eqnarray}\label{term4}
\qquad\frac{\mu_n}{\lambda\sigma_n^2} E_\theta\biggl(K_n \int
_{-K_n/\sigma
_n}^{K_n/\sigma_n} \int_0^1 (1- \mathbf{1}_{F_{n,\theta}})\,du\, dt \biggr)
= \frac{2\mu_n}{\lambda\sigma_n^3} E_\theta[ K_n^2 (1-\mathbf{
1}_{F_{n,\theta}})].
\end{eqnarray}

Substituting the second term in (\ref{indicator2terms}) into (\ref
{befdec}), conditioning on $\mathcal{ F}_n$ and invoking the $\mathcal{ F}_n$
measurability of $K_n$ and $F_{n,\theta}$ provided by Conditions \hyperlink{assumptionKn}{2}
and~\hyperlink{assumptioncondlawYnvgivenK}{3}, respectively, yields
%
%e42 #&#
\begin{eqnarray}\label{int-PJ}
\qquad&&\frac{\mu_n}{\lambda\sigma_n^2} E_\theta\biggl( K_n \int
_{-K_n/\sigma
_n}^{K_n/\sigma_n} \int_0^1 \mathbf{1}(z \le W_n+ut \le z+\lambda)\mathbf{
1}_{F_{n,\theta}}\,du\, dt \biggr)
\nonumber
\\[-8pt]
\\[-8pt]
\nonumber
&&\qquad = \frac{\mu_n}{\lambda\sigma_n^2} E_\theta\biggl( K_n \int
_{-K_n/\sigma_n}^{K_n/\sigma_n} \int_0^1 P_\theta^{\mathcal{ F}_n}(z
\le
W_n+ut \le z+\lambda)\mathbf{1}_{F_{n,\theta}}\,du\, dt \biggr),
\end{eqnarray}
where $P_\theta^{\mathcal{ F}_n}$ denotes conditional probability with
respect to $\mathcal{ F}_n$. To handle the indicator in (\ref{int-PJ}),
note that Condition \hyperlink{assumptioncondlawYnvgivenK}{3} implies that
$n-L_n \ge n_0$ on $F_{n,\theta}$. Hence on $F_{n,\theta}$ we may define
\begin{eqnarray*}
\underline{W_{n,\theta}}=
\frac{V_n-\mu_{n-L_n,\psi_{n,\theta}}}{\sigma_{n-L_n,\psi
_{n,\theta}}}
\end{eqnarray*}
and write
%
%e43 #&#
\begin{eqnarray}\label{Sn=terms}
W_n&=&\biggl( \frac{\sigma_{n-L_n,\psi_n}}{\sigma_n}
\biggr)\underline
{W_n} + \biggl( \frac{Y_n-V_n}{\sigma_n}\biggr)-\biggl( \frac{\mu
_n-\mu
_{n-L_n,\psi_n}}{\sigma_n}\biggr)
\nonumber
\\[-8pt]
\\[-8pt]
\nonumber
&:=& \rho_n \underline{W_n}+T_{n,1}-T_{n,2}.
\end{eqnarray}

By Conditions \hyperlink{assumptionYnmYnVnprime}{5} and \hyperlink{assumptioncondlawYnvgivenK}{3}
we have
$|T_{n,1}| \le B_n/\sigma_n$ and that $\rho_n, B_n$ and $T_{n,2}$ are
$\mathcal{ F}_n$-measurable.
Using (\ref{Sn=terms}) we may write
\begin{eqnarray}\label{PJ}
&& P_\theta^{\mathcal{ F}_n}(z \le W_n + ut \le z+\lambda)\mathbf{
1}_{F_{n,\theta
}} \nonumber\\
&&\qquad=P_\theta^{\mathcal{ F}_n}\bigl( \rho_n^{-1}(z -T_{n,1}+T_{n,2}-ut) \le
\underline{W_n}\nonumber\\
&&\qquad\hspace*{34pt} \le\rho_n^{-1}(z -T_{n,1}+T_{n,2}-ut+\lambda)\bigr)\mathbf{
1}_{F_{n,\theta}}
\nonumber
\\[-8pt]
\\[-8pt]
\nonumber
&&\qquad\le P_\theta^{\mathcal{ F}_n}\bigl( \rho_n^{-1}(z+T_{n,2}-ut)-B_n/\sigma
_{n-L_n,\psi_n}\le\underline{W_n} \\
&&\qquad\hspace*{34pt}\le\rho
_n^{-1}(z+T_{n,2}-ut)+B_n/\sigma_{n-L_n,\psi_n}+\rho_n^{-1} \lambda
\bigr)\mathbf{1}_{F_{n,\theta}} \nonumber\\
&&\qquad= P_\theta^{\mathcal{ F}_n}( Q_n-B_n/\sigma_{n-L_n,\psi_n} \le
\underline
{W_n} \le Q_n+B_n/\sigma_{n-L_n,\psi_n} +\rho_n^{-1}\lambda)\mathbf{
1}_{F_{n,\theta}},\nonumber
\end{eqnarray}
%
%&=& \mathbf{1}( \rho_n^{-1}(z -T_{n,1}+T_{n,2}-ut) \le\underline{W_n}
%&\le& \mathbf{1}( \rho_n^{-1}(z+T_{n,2}-ut)-B_n/\sigma_{n-L_n,\psi_n}\le
%&& \hspace{1cm} \underline{W_n} \le\rho_n^{-1}(z+T_{n,2}-ut)+B_n/
%&=& \mathbf{1}( Q_n-B_n/\sigma_{n-L_n,\psi_n} \le\underline{W_n} \le
%Q_n+B_n/\sigma_{n-L_n,\psi_n} +\rho_n^{-1}\lambda),
where we have set
\[
Q_n= \rho_n^{-1}(z+T_{n,2}-ut).
\]
Recalling (\ref{defdeltanr}), we have
%
%e44 #&#
\begin{eqnarray}\label{eqconcineq}
&&P_\theta(z \le W_{n,\theta} \le z+\lambda)\nonumber\\
&&\qquad\le |P_\theta(z \le W_{n,\theta} \le z+\lambda)-P(z \le Z\le
z+\lambda
)| +P(z \le Z \le z+\lambda)\\
&&\qquad\le 2\delta(n,r_{n,\theta})+\lambda/\sqrt{2 \pi} .\nonumber
\end{eqnarray}
%
%As Condition~\ref{unifvarbnd} holds with $c_2$ and $c_2$ replaced by
%$\max\{c_2,1\}$ and $\max\{c_2,1\}$ %respectively, we may assume that $
% Note that since $\mathcal{ L}_\theta(V_n|J_n)=\mathcal{ L}_{\psi_{n,
%$F_{n,\theta}$ we have
% \bea\label{condcalFn}
% E_\theta[u(\underline{W_n})|\mathcal{ F}_n]=v(n-L_n,\psi_{n,\theta})
% \ena

Since the endpoints of the interval bounding $\underline{W_n}$ in
(\ref
{PJ}) are $\mathcal{ F}_n$-measurable, using Condition~\hyperlink{assumptioncondlawYnvgivenK}{3}
%(\ref{condcalFn}),
and (\ref{eqconcineq}) with the appropriate substitutions, expression
(\ref{PJ}) is bounded by
%(maybe save factor of 2 on $B_V$ by accounting for sign, if $\mathcal{
%J}$-measurable)
%
\begin{eqnarray*}
&&\bigl( 2\delta(n-L_n,r_{n-L_n,\psi_{n,\theta}})+ (2B_n/\sigma
_{n-L_n,\psi_{n,\theta}} + \rho_n^{-1}\lambda)/\sqrt{2 \pi}
\bigr)\mathbf{
1}_{F_{n,\theta}}\\
&&\qquad\le\bigl( 2\delta(n-L_n,r_{n,\theta}/c_2)+\bigl(2\sqrt{c_1} B_n/\sigma
_{n,\theta} + \sqrt{c_1}\lambda\bigr)/\sqrt{2 \pi}\bigr)\mathbf{
1}_{F_{n,\theta}},
\end{eqnarray*}
where we have applied Condition \hyperlink{unifvarbnd}{4},
and that $\delta(n,r)$ is nonincreasing in $r$. As this last quantity
does not depend on $u$ or $t$, substitution into (\ref{int-PJ}) yields
the bound
%
%e45 #&#
\begin{eqnarray}\label{twotermsfromint}
\frac{2\mu_{n,\theta}}{\lambda\sigma_{n,\theta}^3} E_\theta\bigl[
K_n^2 \bigl( 2\delta(n-L_n,r_{n,\theta}/c_2)+ \bigl(2\sqrt{c_1}
B_n/\sigma
_{n,\theta}+ \sqrt{c_1}\lambda\bigr)/\sqrt{2\pi} \bigr)\bigr]\mathbf{
1}_{F_{n,\theta}}.\hspace*{-35pt}
\end{eqnarray}
%
%detailed assumptions about $\sigma_n$.]}
Expression (\ref{twotermsfromint}) leads to three terms. By (\ref
{assumptionKninequality}), that $F_{n,\theta} \subset\{L_n \le
s_{n,\theta}\}$, and since $n-s_{n,\theta} \ge n_0$ for all $\theta
\in
\Theta_{n,r_1}$ implies $s_n=\sup_{\theta\in\Theta
_{n,r_1}}s_{n,\theta
} \le n-n_0$, there exists a positive constant $\Cl{C1.proof}$ such
that the first term satisfies
\begin{eqnarray}
&&\frac{4\mu_{n,\theta}}{\lambda\sigma_{n,\theta}^3} E_\theta[
K_n^2 \delta(n-L_n,r_{n,\theta}/c_2)
\mathbf{1}_{F_{n,\theta}}]\nonumber\\
&&\qquad=\frac{4\mu_{n,\theta}k_{n,\theta,2}}{\lambda\sigma_{n,\theta}^3}
E_\theta\biggl[ \frac{K_n^2}{E_\theta K_n^2} \delta
(n-L_n,r_{n,\theta
}/c_2)\mathbf{1}_{F_{n,\theta}}\biggr]
\nonumber
\\[-8pt]
\\[-8pt]
\nonumber
&&\qquad\le\frac{\Cr{C1.proof}}{\lambda r_{n,\theta}} E_\theta\biggl[
\frac
{K_n^2}{E_\theta K_n^2} \delta(n-L_n,r_{n,\theta}/c_2) \mathbf{
1}_{F_{n,\theta}}\biggr]\\
&&\qquad\le\frac{\Cr{C1.proof}}{\lambda r_{n,\theta}}
\sum_{l=0}^{s_n} \delta(n-l,r_{n,\theta}/c_2)t_{n,\theta,l}, \label
{tn1-recursion}\nonumber
\end{eqnarray}
where $t_{n,\theta,l}$, given in (\ref{deftnthetal}), satisfy
%
%e46 #&#
\begin{equation}\label{snlO1}
\sum_{l=0}^n t_{n,\theta,l} = 1 \qquad\mbox{for all $\theta\in\Theta_{n,r}$.}
\end{equation}

Dropping the indicator $\mathbf{1}_{F_{n,\theta}}$, the sum of the second
and third terms of (\ref{twotermsfromint}) are bounded by
%
%e47 #&#
\begin{equation}\label{tn2-recursion}
\frac{4\sqrt{c_1}\mu_n E_\theta[K_n^2B_n]}{\sqrt{2 \pi} \lambda
\sigma
_n^4} + \frac{\sqrt{2c_1} \mu_n}{ \sqrt{\pi} \sigma_n^3}E_\theta K_n^2.
\end{equation}

Collecting terms (\ref{term1}), (\ref{term2}), (\ref{term4}), (\ref
{tn1-recursion}) and (\ref{tn2-recursion}), and letting
\begin{eqnarray*}
c_{n,\theta,1}&=& \frac{\mu_{n,\theta}}{\sigma_{n,\theta}^2}\Psi
_{n,\theta} +
\frac{\mu_{n,\theta}}{\sigma_{n,\theta}^3}
E_\theta\biggl[ \biggl( \biggl( 1 +\frac{\sqrt{2c_1}}{\sqrt{\pi}}
\biggr)+|W_{n,\theta}| \biggr)K_n^2\biggr]\quad
\mbox{and}\\
c_{n,\theta,2}&=& \frac{2\mu_{n,\theta}}{\sigma_{n,\theta}^3}
E_\theta
[ K_n^2(1-\mathbf{1}_{F_{n,\theta}})]+\frac{4\sqrt{c_1}\mu
_{n,\theta} E_\theta[K_n^2B_n]}{\sqrt{2 \pi} \sigma_{n,\theta}^4},
\end{eqnarray*}
for all $z \in\mathbb{R}$ we have
%
%e48 #&#
\begin{eqnarray}\label{Ehzlambda}
&&|E_\theta h_{z,\lambda}(W_{n,\theta})-Nh_{z,\lambda}|
\nonumber
\\[-8pt]
\\[-8pt]
\nonumber
&&\qquad \le\frac{\Cr
{C1.proof}}{\lambda r_{n,\theta}}\sum_{l=0}^{s_n} \delta
(n-l,r_{n,\theta
}/c_2)t_{n,\theta,l} + c_{n,\theta,1}+ \frac{1}{\lambda}c_{n,\theta,2}.
\end{eqnarray}

Note that Conditions \hyperlink{Psinisrootn}{1} and \hyperlink{assumptionKn}{2}, and
\hyperlink{assumptioncondlawYnvgivenK}{3}
and \hyperlink{assumptionYnmYnVnprime}{5}, respectively, yield the existence of
positive constants $\Cl{C2.proof}$ and $\Cl{C3.proof}$ that
%
%e49 #&#
\begin{eqnarray}\label{sees}
c_{n,\theta,1} \le\Cr{C2.proof}/r_{n,\theta} \quad\mbox{and}\quad
c_{n,\theta,2}
\le\Cr{C3.proof}/r_{n,\theta}^2.
\end{eqnarray}
As $\mathbf{1}(w \le z) \le h_{z,\lambda}(w) \le\mathbf{1}(w \le z+\lambda)$
we obtain
\begin{eqnarray*}
&&P_\theta(W_{n,\theta} \le z) - P(Z \le z) \\
&&\qquad\le |E_\theta
h_{z,\lambda
}(W_{n,\theta}) - Eh_{z,\lambda}(Z)| + Eh_{z,\lambda}(Z)-P(Z \le z)
\end{eqnarray*}
with $Eh_{z,\lambda}(Z)-P(Z \le z) \le P(z \le Z \le z+\lambda) \le
\lambda/\sqrt{2 \pi}$. Along with a similar lower bound obtained by
considering $h_{z-\lambda,\lambda}(w)$, in view of (\ref{Ehzlambda})
and (\ref{sees}) we have that for every $z \in\mathbb{R}$
\begin{eqnarray*}
&&|P_\theta(W_{n,\theta} \le z) - P(Z \le z)| \\
&&\qquad\le\frac{\Cr
{C1.proof}}{\lambda r_{n,\theta}}\sum_{l=0}^{s_n} \delta
(n-l,r_{n,\theta
}/c_2)t_{n,\theta,l} + \frac{\Cr{C2.proof}}{r_{n,\theta}}+
\frac{\Cr{C3.proof}}{\lambda r_{n,\theta}^2} + \frac{\lambda
}{\sqrt{2
\pi}}.
\end{eqnarray*}

Letting $\lambda= 2c_2\Cr{C1.proof}/r_{n,\theta}$, and, noting that
the right-hand side does not depend on~$z$, taking supremum over $z \in
\mathbb{R}$ yields
\begin{eqnarray}\label{bdinthetarntheta}
&&\sup_{z \in\mathbb{R}}|P_\theta(W_{n,\theta}\le z)-P(Z \le z)|\nonumber\\
&&\qquad\le\sum_{l=0}^{s_n} \delta(n-l,r_{n,\theta}/c_2)t_{n,\theta
,l}/2c_2 +
\Cl{b.proof}/r_{n,\theta}\\
&&\qquad\le\sum_{l=0}^{s_n} \delta(n-l,r/c_2)t_{n,\theta,l}/2c_2 +
\Cr{b.proof}/r\nonumber
\end{eqnarray}
for $\Cr{b.proof} = \Cr{C2.proof} + \Cr{C3.proof}/2c_2\Cr{C1.proof}
+2c_2\Cr{C1.proof}/\sqrt{2 \pi}$, where for the last inequality we have
used that $\theta\in\Theta_{n,r}$, and that $\delta(n,r)$ and $1/r$
are nonincreasing functions of $r$. Taking supremum over $\Theta
_{n,r_1}$ on the right-hand side of (\ref{bdinthetarntheta}), then
over $\Theta_{n,r} \subset\Theta_{n,r_1}$ on the left yields
%
%e50 #&#
\begin{equation}\label{bdinthetarnthetal}
\delta(n,r) \le
\sup_{\theta\in\Theta_{n,r_1}} \sum_{l=0}^{s_n} \delta
(n-l,r/c_2)t_{n,\theta,l}/2c_2 + \Cr{b.proof}/r.
\end{equation}

Suppose first that Condition \hyperlink{Thetanr-compact}{6}a is satisfied, so
that $L_n=l_{n,0}$ almost surely for some $l_{0,n} \in\{0,\ldots,n\}$
for all $\theta\in\Theta_{n,r_1}$. If $l_{0,n}>s_n$ then (\ref
{deftnthetal}) and (\ref{bdinthetarnthetal}) yield $\delta(n,r)
\le\Cr{b.proof}/r$, proving (\ref{pfgoal}). Otherwise $t_{n,\theta
,l}=\mathbf{1}(l=l_{n,0})$ for $0 \le l_{n,0} \le s_n$, and inequality
(\ref{bdinthetarnthetal}) specializes to
%
%e51 #&#
\begin{equation}\label{bdinthetarnthetal0}
\delta(n,r) \le
\delta(n-l_{n,0},r/c_2)/2c_2 + \Cr{b.proof}/r.
\end{equation}

When Condition \hyperlink{Thetanr-compact}{6}b is satisfied, the sum in (\ref
{bdinthetarnthetal}) is a continuous function of~$\theta$ on the
compact set $\Theta_{n,r_1}$, and hence achieves its supremum at some
$\theta_n^* \in\Theta_{n,r_1}$. Letting $p_{n,l}=t_{n,\theta
_n^*,l}/2$, from (\ref{bdinthetarnthetal}) and (\ref{snlO1}) we have
%
%e52 #&#
\begin{eqnarray}\label{proofdelta0nr0}
\delta(n,r) \le\sum_{l=0}^{s_n} \delta(n-l,r/c_2)p_{n,l}/c_2 + \Cr
{b.proof}/r \qquad\mbox{with } \sum_{l=0}^n p_{n,l}=1/2.
\end{eqnarray}
As (\ref{bdinthetarnthetal0}) is the special case of (\ref
{proofdelta0nr0}) when $p_{n,l}=\mathbf{1}(l=l_{n,0})/2$, it suffices to
handle the latter.

Let $a_n=0$ for $0 \le n <n_0$, and $a_n=\sup_{r \ge r_1} r\delta(n,r)$
for $n \ge n_0$. For all $r \ge r_1$ and $n \ge n_0$ we have
\begin{eqnarray*}
(r/c_2)\delta(n,r/c_2) &\le& \sup_{s\dvtx  s \ge r_1} (s/c_2)\delta(n,s/c_2)
\\
&=& \sup_{s\dvtx  s \ge r_1/c_2} s\delta(n,s) \\
&\le& \Bigl[\sup_{s\dvtx  r_1/c_2 \le s < r_1} s\delta(n,s)\Bigr]\mathbf{
1}(c_2>1) + \sup_{s\dvtx s \ge r_1} s\delta(n,s)\\
&\le& r_1+a_n.
\end{eqnarray*}

Using that $n \ge n_1$ implies $n-s_n \ge n_0$, multiplication by $r$
in (\ref{proofdelta0nr0}) yields, with $f=r_1/2+\Cr{b.proof}$, that
for all $n \ge n_1$
\begin{eqnarray*}
r\delta(n,r) &\le&\sum_{l=0}^{s_n} (r/c_2)\delta(n-l,r/c_2)p_{n,l} +
\Cr{b.proof}
\le\sum_{l=0}^{s_n}(r_1+a_{n-l})p_{n,l} + \Cr{b.proof} \\
&\le&
\sum_{l=0}^n a_{n-l}p_{n,l}+f.
\end{eqnarray*}
Taking supremum on the left-hand side over $r \ge r_1$ and recalling
(\ref{proofdelta0nr0}) now yields
\begin{eqnarray*}
a_n \le
\sum_{l=0}^n a_{n-l} p_{n,l} + f  \qquad\mbox{with }  \sum_{l=0}^n
p_{n,l}=1/2  \mbox{ for all $n \ge n_1$.}
\end{eqnarray*}
Lemma~\ref{big-recursion} now implies $\sup_{n \ge n_1} a_n < \infty$.
Hence, there exists a constant $C$ such that
$\delta(n,r) \le C/r$ for all $n \ge n_1$ and all $r \ge r_1$; that is,
(\ref{pfgoal}) holds.
\end{pf*}

%

%in a combinatorial central limit theorem. \textit{Z.
%Wahrscheinlichkeitstheorie verw. Gebiet.}, {\bf66}, pp 379-386.

%combinatorial central limit theorems and pattern occurrences,
%using zero and size biasing. \textit{Jour. Appl. Probab.}, {\bf42},
%pp. 661-683. arXiv:math.PR/0511510

%approximation for coverage
% models over binomial point processes

% \item\label{GR} Goldstein, L, and Rinott, Y. (1996)
% Multivariate normal approximations by Stein's method and size bias
%couplings,

%approximation with Stein's method of exchangeable pairs under a
%general linearity condition,
%arXiv:0711.1082v2 [math.PR]

%constructions and rates in the CLT for dependent
%summands with applications to the antivoter model and weighted
%$U$-statistics. \textit{Ann. Appl. Probab.}, {\bf7}, 1080-1105.

%
%expectations.

% imsref loaded by akundreckaite, 2012-05-16 10:05:17
%

%suskaldyti doi

\printaddresses

\end{document}